\documentclass[11pt]{article}
\usepackage[T1]{fontenc}
\usepackage{lmodern}
\usepackage[margin=1in]{geometry}
\usepackage{amsmath,amssymb,amsthm,mathtools}
\usepackage{microtype}
\usepackage{xcolor}
\usepackage{tikz}
\usepackage{subfig}
\usepackage[hidelinks]{hyperref}
\usepackage{xcolor}
\numberwithin{equation}{section}
\newtheorem{theorem}{Theorem}[section]
\theoremstyle{definition}

\newtheorem{proposition}[theorem]{Proposition}
\newtheorem{lemma}[theorem]{Lemma}

\theoremstyle{remark}
\newtheorem{remark}[theorem]{Remark}
\newcommand{\E}{\mathcal E}
\newcommand{\F}{\mathcal F}
\newcommand{\cL}{\mathcal L}

\newcommand{\PP}{\mathbb P}
\newcommand{\EE}{\mathbb E}
\newcommand{\ind}{\mathbf 1}
\newcommand{\dd}{\,\mathrm d}

\newcommand{\muae}{\text{a.e.}}
\newcommand{\qe}{\text{q.e.}}

\newcommand{\R}{\mathbb R}
\newcommand{\N}{\mathbb N}
\newcommand{\cE}{\mathcal E}
\newcommand{\cF}{\mathcal F}

\renewcommand{\thefootnote}{}

\AtEndDocument{
	\bigskip{\footnotesize
		\textsc{School of Mathematics, Hangzhou Normal University, Hangzhou 310036, China} \par \textsc{School of Mathematics, Nanjing University, Nanjing 210093, China.}\par
		\textit{E-mail address}: \texttt{gaojin@hznu.edu.cn} \par
	}
	\bigskip{\footnotesize%
	\textsc{Department of Mathematical Sciences, Tsinghua University, Beijing 100084, China}  \par  
	\textit{E-mail address}: \texttt{songyj24@mails.tsinghua.edu.cn} \par
	}
	\bigskip{\footnotesize
		\textsc{College of Science, National University of Defense Technology, Changsha 410073, China}  \par
		\textit{E-mail address}: \texttt{yuzy23@nudt.edu.cn} \par
	}
	\bigskip{\footnotesize
		\textsc{School of Mathematics, South China University of Technology, Guangzhou 510641, China.}  \par
		\textit{E-mail address}: \texttt{summerfish@scut.edu.cn} \par
	}

}

\title{Some open problems on sub-Gaussian upper estimates of heat kernels}

\author{Jin Gao, Yijun Song, Zhenyu Yu, and Junda Zhang}
\date{}

\newcommand{\Capm}{\operatorname{cap}}

\newcommand{\HK}{\operatorname{HK}}

\newcommand{\FK}{\operatorname{FK}}

\newcommand{\VC}{\operatorname{V}}

\begin{document}
\maketitle

\begin{abstract}
We present answers to two open problems on sub-Gaussian upper estimates of heat kernels on metric measure spaces.

The first, posed by Grigor'yan, Hu and Lau, conjectures that a near-diagonal lower estimate implies the global sub-Gaussian upper estimate under a lower Ahlfors regularity condition. We prove it by introducing an annular resolvent estimate.

The second, also due to Grigor'yan, Hu and Lau, conjectures that the generalized capacity condition in the characterization of the heat kernel upper estimate can be replaced by a capacity upper bound. We present a concrete counterexample: the relative Faber--Krahn inequality and the capacity upper bound hold, whereas the cutoff Sobolev inequality and the sub-Gaussian upper estimate fail.
\end{abstract}

\footnote{\textsl{MSC2020}: 28A80, 35K08}
\footnote{\textsl{Keywords}: Heat kernel, sub-Gaussian upper estimates, capacity upper bound, Faber--Krahn inequality.}
\renewcommand{\thefootnote}{\arabic{footnote}}
\setcounter{footnote}{0}


\section{Introduction}

The study of heat kernel estimates is a central topic in analysis on metric measure spaces. Let $(X,d,m)$ be a metric measure space, and let $\{P_t\}_{t>0}$ be the semigroup associated with a regular Dirichlet form $(\mathcal E,\mathcal F)$ on $L^2(X,m)$. Denote by $p_t(x,y)$ the corresponding heat kernel. A fundamental object is the sub-Gaussian upper estimate
\[
p_t(x,y)\le \frac{C}{m(B(x,t^{1/\beta}))}
\exp\left(-c\left(\frac{d(x,y)^\beta}{t}\right)^{1/(\beta-1)}\right)
\]
for $t>0$ and $m$-a.e. $x,y\in X$, where $\beta>1$. Such estimate is intimately connected with parabolic Harnack inequalities, Poincar\'e inequalities, cutoff Sobolev inequalities, Faber--Krahn inequalities, and capacity conditions. This paper is concerned with the following related open problems.

\begin{enumerate}
	\item \textbf{Conjecture 1} (\cite[Conjecture 3.6]{GHL08}). Under a lower Ahlfors regular measure condition, the near-diagonal lower estimate of the heat kernel implies the global sub-Gaussian upper estimate.

	\item \textbf{Conjecture 2} (Generalized capacity condition in \cite[(1.36)]{GHL15} can be replaced by the capacity upper bound). Let $(M,d,\mu)$ be a metric measure space with precompact metric balls. Let $(\E,\F)$ be a regular, strongly local Dirichlet form in $L^2(M,\mu)$ and $\Psi$ be a scaling function. If condition (VD) is satisfied, then the following holds
	\begin{equation*}
		\mathrm{FK}(\Psi)+\mathrm{cap}_{\leq}(\Psi)\implies \mathrm{Gcap}_{\leq} (\Psi)
		\end{equation*}
\end{enumerate}

We provide a positive answer to Conjecture 1 and a concrete counterexample to Conjecture 2. We roughly describe our strategy for these problems.

For Conjecture 1, the key idea is a new annular estimate, which is the distinction from \cite[Section~4]{GHL08}. There, the passage to the full upper bound uses a near-diagonal lower bound for the heat kernel estimate on a ball. Here, each annulus contributes a uniform amount to certain occupation integral after its first crossing, and summing these contributions yields the desired estimate. This can be regarded as counting annular passages through their occupation contributions, including possible returns.

For Conjecture 2, we construct an unbounded complete geodesic Ahlfors $2\alpha$-regular metric measure space $(X,d,m)$ (based on the unbounded Vicsek fractal) and a strongly local regular Dirichlet form $(\mathcal E,\mathcal F)$ on $L^2(X,m)$ with the following properties. The relative Faber--Krahn inequality $\mathrm{FK}(\beta)$ and the capacity upper bound $\mathrm{cap}_{\leq}(\beta)$ hold, whereas the generalized capacity inequality
$\mathrm{Gcap}_{\leq}$ and the sub-Gaussian upper estimate fail. Thus the generalized capacity condition in the upper-bound characterization of \cite[Theorem~1.3]{GHL15} cannot be replaced by the capacity upper bound.

The basic idea of the construction is to glue two spaces with different walk-dimensions in order to fail the Poincar\'e inequality. We use the Vicsek set simply because it is easy to analyze among classical fractals. Let $Y$ be the unbounded Vicsek set and let $Z$ be the chain of unit Vicsek cells along a diagonal ray. We take
\[
S=Y\times Y,\qquad F=Z\times Z,
\]
and glue them along the countably many diagonal vertices. The direct products increase the Hausdorff dimension from $\alpha$ to $2\alpha$ without changing the small-scale walk-dimension $\beta$. Since $2\alpha>\beta$, diagonal vertices have zero $1$-capacity and so their values do not affect the Dirichlet-form energy. This is the main purpose of taking products.

The paper is organized as follows. In Section~2, we state the precise definitions and main results. Sections~3 and~4 are devoted to the proofs of the two main theorems, respectively.

\section{Setting and main results}

	Let $(M,d)$ be a nonempty locally compact separable metric space, and let
	$\mu$ be a Radon measure with full support. Let $(\E,\F)$ be a local regular
	symmetric Dirichlet form on $L^2(M,\mu)$. Here locality means that the form has no jump part, it does not by itself exclude a killing measure. No conservativeness assumption is imposed, as in the setting of \cite[Conjecture~3.6]{GHL08}.
	Write
	\[
	B(x,r)=\{y\in M:d(x,y)<r\},\qquad V(x,r)=\mu(B(x,r)).
	\]
	In particular, $\mu$ has full support and is $\sigma$-finite. Let $\cL$
	be its nonnegative self-adjoint generator on the Hilbert space
	$L^2(M,\mu)$, and let $P_t=e^{-t\cL}$ be its corresponding semigroup
	 for $t\ge0$. We write $\cE_1(u,v)=\cE(u,v)+(u,v)_{L^2(M,\mu)}$.  For an open set $U\subset M$, let
	\begin{equation}\label{eq:F0U}
		\cF(U)=\overline{\cF\cap C_c(U)}^{\,\cE_1}.
	\end{equation}

	Fix $\alpha>0$, $\beta>1$, and $r_0\in(0,\infty]$. We say that the upper $\alpha$-regularity condition \eqref{eq:upper-volume} holds if there exists
	a constant $C_1$ such that, for $\mu$-almost all $x\in M$ and all $0<r<r_0$,
	\begin{equation}\label{eq:upper-volume}
		\mu(B(x,r))\leq C_1 r^\alpha. \tag*{V$_{\leq}$}
	\end{equation}

  We say that the lower $\alpha$-regularity condition \eqref{eq:lower-volume} holds if there exists
a constant $c_V>0$ such that, for $\mu$-almost all $x\in M$ and all $0<r<r_0$,
	\begin{equation}\label{eq:lower-volume}
		\mu(B(x,r))\ge c_V r^\alpha,
		\tag*{V$_{\geq}$}
	\end{equation}
  We say that $\mu$ is $\alpha$-regular ($V_{\alpha}$) if $\mu$ is both
  upper and lower $\alpha$-regular.

 We say that the relative Faber--Krahn inequality \eqref{eq:FK} holds if there are $C_F>0$ and $\nu>0$
such that, for every ball $B=B(x,r)$ and every nonempty open
$U\subset B$,
\begin{equation}\label{eq:FK}
	\lambda_{\min}(U)\geq
	\frac{C_F}{r^\beta}
	\left(\frac{\mu(B)}{\mu(U)}\right)^\nu,\tag*{FK}
\end{equation}
where the bottom of the Dirichlet spectrum of $U$ is
\[
\lambda_{\min}(U):=\inf\left\{\frac{\E(f,f)}{\|f\|_2^2}: f\in\F(U)\setminus\{0\}\right\}.
\]

For any subset
$A\subseteq X$, its $1$-capacity is defined by
$$
\Capm_1(A):=\inf\left\{\cE_1(u,u): u\in\cF,\ \tilde u\ge 1 \text{ q.e. on } A\right\}.
$$
where $\tilde u$ is a quasi-continuous modification of $u$. 
For Borel sets $A_1,A_2\subseteq X$, define
\[
\Capm(A_1,A_2):=\inf\left\{\cE(\varphi,\varphi): \varphi\in\cF,\
\begin{array}{l}
	\varphi=1 \text{ in an open neighborhood of } A_1,\\
	\varphi=0 \text{ in an open neighborhood of } A_2
\end{array}\right\}.
\]

We say that the two-sided capacity bounds $(\Capm)$ hold if both the capacity upper
bound $(\Capm_{\leq})$ and the capacity lower bound $(\Capm_{\geq})$ hold as
follows. There exist $C>0$ and $A>1$ such that for any ball $B(x,r)$,
\begin{equation*}
	\Capm\left(B(x,r),B(x,Ar)^c\right)
	\le C\frac{\mu(B(x,r))}{r^{\beta}},
	\tag*{cap$_{\leq}$}
\end{equation*}
\begin{equation*}
	\Capm\left(B(x,r),B(x,Ar)^c\right)
	\ge\frac{1}{C}\frac{\mu(B(x,r))}{r^{\beta}}.
	\tag*{cap$_{\geq}$}
\end{equation*}


		We say that generalized capacity inequality
		$(\mathrm{Gcap}_{\leq})$ hold
		 if there exist
		constants $C_1,C_2>0$ such that, for any
		$f\in\cF\cap L^\infty(X)$ and any two concentric
		balls
		\[
		B_1=B(x,R),\qquad B_2=B(x,R+r),
		\qquad R,r>0,
		\]
		there exists a cutoff function $\phi\in\cF$
		for $B_1\subseteq B_2$ satisfying
		\begin{align}
		\int_{B_2\setminus B_1}
		|\tilde f|^2\,d\Gamma(\phi,\phi)
		\leq
		C_1\int_{B_2\setminus B_1}d\Gamma(f,f)
		+\frac{C_2}{r^\beta}
		\int_{B_2\setminus B_1}|f|^2\,dm.\label{gcap}\tag*{$\mathrm{Gcap}_{\leq}$}
		\end{align}
		Here $0\leq\phi\leq1$, $\phi=1$ on an open
		neighborhood of $B_1$, and
		$\operatorname{supp}\phi\Subset B_2$.
		The cutoff function $\phi$ may depend on $f$,
		whereas $C_1,C_2$ are independent of
		$x,R,r$ and $f$.

Here $\Gamma$ denotes the energy measure of the strongly local form, normalized by $\mathcal E(f,f)=\int_M d\Gamma(f)$. The energy measure of a cutoff vanishes outside its support. See \cite{AB15,BBK06} for related cutoff Sobolev inequalities.

 We say that the on-diagonal upper estimate of the heat kernel \eqref{eq:DUE} holds if there exists a constant $C>0$ such that  for $0<t<r_0^\beta$ and for $\mu\times\mu$-almost every $x,y\in M$, we have
\begin{equation}\label{eq:DUE}
	p_t(x,y)\leq C t^{-\alpha/\beta}.\tag*{DUE}
\end{equation}

We say that the upper estimate of the heat kernel \eqref{eq:UE} holds if there are constants $C,c>0$ such that for $0<t<r_0^\beta$ and for $\mu\times\mu$-almost every $x,y\in M$, we have
\begin{equation}\label{eq:UE}
	p_t(x,y)\le \frac{C}{t^{\alpha/\beta}}
	\exp\!\left(-c\left(\frac{d(x,y)}{t^{1/\beta}}\right)^{\frac{\beta}{\beta-1}}\right),
	\tag*{UE}
\end{equation}
In particular, when $r_0=\infty$, the estimate holds for every $t>0$.

 We say that the near-diagonal lower estimate of the heat kernel \eqref{eq:NLE} holds if there exist constants $\delta\in(0,1)$ and $c_L>0$ such that for $0<t<\delta r_0^\beta$ and  for $\mu\times\mu$-almost every $x,y\in M$ with $d(x,y)<\delta t^{1/\beta}$, we have
	\begin{equation}\label{eq:NLE}
		p_t(x,y)\ge c_L t^{-\alpha/\beta}.\tag*{NLE}
	\end{equation}

	The following theorem verifies Conjecture 1 (see \cite[Conjecture 3.6]{GHL08}).
	\begin{theorem}\label{thm:main1}
		Let $(M,d,\mu)$ be a metric measure space satisfying \eqref{eq:lower-volume}.
		Assume that the heat kernel $p_t(x,y)$ of the Dirichlet form $(\E,\F)$ exists and satisfies \ref{eq:NLE}. Then \ref{eq:UE} holds.
	\end{theorem}

The near-diagonal lower estimate already yields the on-diagonal upper estimate by \cite[Theorem~3.3]{GHL08}. We use probabilistic arguments to pass from this estimate to the global upper bound.

The following theorem disproves Conjecture 2.

In the counterexample we write $(X,d,m)$ instead of $(M,d,\mu)$. We write $\VC(\Phi)$ if $m(B(x,r))\asymp\Phi(r)$ uniformly in $x$ and $r>0$, and $\VC(q)$ if $\Phi(r)=r^q$. A cutoff for $B(x,r)\subset B(x,Ar)$ is a function in the form domain with values in $[0,1]$, equal to $1$ quasi-everywhere on the inner ball and supported in the outer ball. For an increasing scale function $\Psi$, the notions $\FK(\Psi)$, $\Capm(\Psi)_{\leq}$, $\Capm(\Psi)_{\geq}$, $\Capm(\Psi)$ mean that $r^\beta$ in the preceding definitions are replaced by $\Psi(r)$. We write $\FK(\beta)$, $\Capm(\beta)_{\leq}$, $\Capm(\beta)_{\geq}$, $\Capm(\beta)$ for the power scale $\Psi(r)=r^\beta$, and similarly for the capacity and heat kernel conditions. We use $\HK(\Psi)$ for the corresponding two-sided heat kernel estimates, with prefactor $m(B(x,\Psi^{-1}(t)))^{-1}$ and exponential term
\[
\exp\bigl(-c\Upsilon(Cd(x,y),t)\bigr),\qquad
\Upsilon(R,t)=\sup_{s>0}\left(\frac{R}{s}-\frac{t}{\Psi(s)}\right),
\]
allowing different positive constants in the upper and lower bounds. These conventions agree with the standard notation when $\Psi(r)=r^\beta$. Constants implicit in $\lesssim$ and $\asymp$ are independent of the displayed variables. Write $\N=\{0,1,2,\ldots\}$ and $l(A)$ for all real-valued functions on $A$.

\begin{theorem}\label{thm:main2}
There exist an unbounded complete geodesic metric measure space $(X,d,m)$ and a
strongly local regular Dirichlet form $(\cE,\cF)$ on $L^2(X;m)$ such that
$\VC(2\alpha)$, $\FK(\beta)$, and $\Capm(\beta)_{\leq}$ hold, while $\mathrm{Gcap}_{\leq}$ fails, where
\begin{equation}\label{albe}
 \alpha=\frac{\log 5}{\log 3},\qquad
\beta=\alpha+1=\frac{\log 15}{\log 3}.
\end{equation}
\end{theorem}

	\section{Proof of Theorem \ref{thm:main1}}
	\label{sec:probability}

	Let $X=(X_t,\PP_x)$ be a Hunt process properly associated with
	$(\E,\F)$, with lifetime $\zeta$ (see \cite[Theorems~7.2.1--7.2.2]{FOT}).

	The transition semigroup is denoted by
	\[
	P_tf(x)=\EE_x[f(X_t);\ t<\zeta].
	\]
	Locality of the Dirichlet form implies that, outside
	a properly exceptional set, the paths are continuous on $[0,\zeta)$.
	For $\lambda>0$ and a bounded nonnegative $f\in L^2$, the analytic
	resolvent
	\[
	G_\lambda f=(\cL+\lambda)^{-1}f\in\F
	\]
	has the quasi-continuous representative
	\begin{equation}\label{eq:resolvent-representation}
		G_\lambda f(x)=\EE_x\int_0^\zeta
		e^{-\lambda s}f(X_s)\dd s,
		\qquad \qe\ x.
	\end{equation}
	In particular, $0\le G_\lambda f\le\lambda^{-1}\|f\|_\infty$.
	Every capacity-zero exceptional set used below is included in a
	properly exceptional set. Starting outside such a set, the process
	almost surely avoids it at all times before the lifetime.

	The following lemma immediately follows from the uniqueness of quasi-continuous representatives (see \cite[Section~2.1]{FOT}). Capacity and quasi-continuity are taken with respect to $(\E,\F)$.

	\begin{lemma}\label{lem:local-order}
		Let $U\subset M$ be open, $u\in\F\cap L^\infty$, and $b\in\mathbb R$.
		If $u\ge b$ $\mu$-almost everywhere on $U$, then the quasi-continuous
		representative $\widetilde u$ satisfies $\widetilde u\ge b$ quasi-everywhere
		on $U$.
	\end{lemma}

\begin{proof}

Since polar sets are $\mu$-null, so
$\widetilde u=u$ $\mu$-a.e.\ on $M$. Hence
$g:=\widetilde u-b$ is quasi-continuous on $M$ and $g\ge0$
$\mu$-a.e.\ on $U$. By \cite[Lemma~2.1.4]{FOT}, $g\ge0$ q.e.
on $U$, that is, $\widetilde u\ge b$ q.e. on $U$.
\end{proof}

Without loss of generality, decrease the constant in \ref{eq:NLE} so that $0<\delta\le1/2$. Let
\[
\varepsilon=\frac18\delta^{1+1/\beta},\qquad R_*=\varepsilon r_0,
\]
with $R_*=\infty$ when $r_0=\infty$. In \cite[proof of Theorem~3.2]{GHL08}, it is proved that \ref{eq:NLE} implies
\begin{equation}\label{eq:local-upper-volume}
\mu(B(x_0,R))\le C R^\alpha,\qquad x_0\in M,\quad 0<R<R_*.
\end{equation}

	For $x_0\in M$ and $R>0$, define the stopping times 
	\begin{equation}\label{eq:crossing-time}
		S(x_0,R)=\inf\{s<\zeta:d(x_0,X_s)\ge R\},
		\qquad \inf\varnothing=\infty.
	\end{equation}
	If the process starts in $B(x_0,R)$ with $S(x_0,R)<\infty$, then the continuity gives
	\begin{equation}\label{eq:crossing-location}
		S(x_0,R)<\zeta,
		\qquad d(x_0,X_{S(x_0,R)})=R.
	\end{equation}
	We use the convention $e^{-\lambda\infty}=0$ throughout this paper.

	We first establish the following important annular resolvent estimate.

	\begin{lemma}\label{lem:fixed-ball}
		There exists a constant $a>0$, such that, for each $x_0\in M$ and $0<R<R_*$, one can find a properly
		exceptional set $N_{x_0,R}$ satisfying
		\begin{equation}\label{eq:fixed-ball-laplace}
			\EE_y\exp\!\left[-\frac{S(x_0,R)}{aR^\beta}\right]\le\frac12,
			\qquad y\in B(x_0,R/4)\setminus N_{x_0,R}.
		\end{equation}
	\end{lemma}

	\begin{proof}
		Let
		\begin{equation}\label{eq:q-and-b}
			q=c_Lc_V\delta^\alpha,
			\qquad b=q(e^{-1}-e^{-2})>0,
		\end{equation}
		and choose a large integer
		$N\ge\max\{1,2/b\}$. Fix $x_0$ and $R$. Let
		\[
		h=\frac{R}{8N},
		\qquad \rho_i=\frac R2+(2i-1)h,
		\qquad i=1,\ldots,N,
		\]
		and define
		\begin{align*}
			A_i&=\{z\in M:|d(x_0,z)-\rho_i|<h/2\},\\
			U_i&=\{z\in M:|d(x_0,z)-\rho_i|<h/4\}.
		\end{align*}
		The sets $A_i$ are pairwise disjoint and contained in $B(x_0,R)$. Let
		\begin{equation}\label{eq:T-and-a}
			T=(h/8)^\beta=aR^\beta,
			\qquad a=(64N)^{-\beta},
			\qquad \lambda=T^{-1}.
		\end{equation}
		For $z\in U_i$ and $T\le s\le2T$, we have
		\begin{equation}\label{eq:ball-in-annulus}
			B(z,\delta s^{1/\beta})\subset A_i,
		\end{equation}
		since $\delta s^{1/\beta}\le\delta 2^{1/\beta}h/8<h/4$ and also $2T<\delta r_0^\beta$, by the definition of $R_*$.

		By \eqref{eq:local-upper-volume}, the functions $\ind_{A_i}$ belong to
		$L^2$. Let
		\[
		u_i=G_{1/T}\ind_{A_i}\in\F\cap L^\infty.
		\]
		By \ref{eq:NLE}, \eqref{eq:ball-in-annulus}, and
		\ref{eq:lower-volume},
		\[
		P_s\ind_{A_i}(z)
		\ge c_Ls^{-\alpha/\beta}
		\mu(B(z,\delta s^{1/\beta}))
		\ge q
		\]
		for almost every $z\in U_i$ and each $s\in[T,2T]$.
		Therefore
		\begin{equation}\label{eq:annular-resolvent-lower}
			u_i(z)\ge q\int_T^{2T}e^{-s/T}\dd s=bT,
			\qquad \mu\text{-}\muae\ z\in U_i.
		\end{equation}
		Lemma~\ref{lem:local-order} gives the same bound quasi-everywhere
		on $U_i$. Enlarge a properly exceptional set, so that the
		probabilistic representations of all the $u_i$, the quasi-everywhere lower bounds,
		and path continuity remain valid off that set.

		Start at $y\in B(x_0,R/4)$ outside this enlarged exceptional set. Let $\sigma_i$ be the first
		hitting time before $\zeta$ of the sphere $\{d(x_0,\cdot)=\rho_i\}$, with
		$\sigma_i=\infty$ when no hit occurs. Let $\sigma=S(x_0,R)$.
		On $\{\sigma<\infty\}$, path continuity gives
		$\sigma_i\le\sigma<\zeta$ for every $i$. On
		$\{\sigma_i<\infty\}$ the hit point lies in $U_i$ and outside the
		enlarged exceptional set. The strong Markov property and the
		nonnegativity of the occupation before $\sigma_i$ imply
		\begin{align}
			u_i(y)
			\ge\EE_y\!\left[
			e^{-\sigma_i/T}u_i(X_{\sigma_i});\ \sigma_i<\zeta
			\right]\ge bT\,\EE_y[e^{-\sigma_i/T}]
			\ge bT\,\EE_y[e^{-\sigma/T}].
			\label{eq:resolvent-stopping}
		\end{align}
		All expressions on paths with an infinite crossing time are understood
		to be zero.

		On the other hand, using the probabilistic representations and the
		disjointness of the annuli,
		\begin{align}
			\sum_{i=1}^N u_i(y)=\EE_y\int_0^\zeta e^{-s/T}
			\sum_{i=1}^N\ind_{A_i}(X_s)\dd s\le\int_0^\infty e^{-s/T}\dd s=T.
			\label{eq:sum-resolvents}
		\end{align}
		Adding \eqref{eq:resolvent-stopping} and applying
		\eqref{eq:sum-resolvents} gives
		\[
		\EE_y[e^{-\sigma/T}]\le\frac1{Nb}\le\frac12,
		\]
		proving the assertion.
	\end{proof}

	\begin{remark}\label{rem:core}
		Passing from almost everywhere to quasi-everywhere in
		\eqref{eq:annular-resolvent-lower} is essential, since the sphere-hitting
		distributions are not known to be absolutely continuous with respect to $\mu$.

	\end{remark}

	As in \cite[Section~4]{GHL08}, we construct a common exceptional set valid for all centers and radii.

	\begin{proposition}\label{prop:uniform-laplace}
		There are a properly exceptional set $\mathcal N$ and a constant $K>0$
		such that
		\begin{equation}\label{eq:uniform-laplace}
			\EE_x[e^{-Kr^{-\beta}S(x,r)}]\le\frac12,
			\qquad x\in M\setminus\mathcal N,\quad 0<r<R_*.
		\end{equation}
	\end{proposition}

	\begin{proof}
		Take a countable dense subset $D\subset M$ and apply
		Lemma~\ref{lem:fixed-ball} to the countable family
		\[
		\{B(x_0,R):x_0\in D,\ R\in\mathbb Q\cap(0,R_*)\}.
		\]
		We put the union of the corresponding exceptional sets, and the set
		needed for path continuity, into one properly exceptional set
		$\mathcal N$. This is possible because a countable union of polar sets is polar and is contained in a properly exceptional set; see \cite[Chapter~4]{FOT}. In particular,
\[
\PP_x\bigl(X_s\in\mathcal N\text{ for some }s<\zeta\bigr)=0,
\qquad x\notin\mathcal N.
\]
Thus the estimates remain applicable after the stopping times used below.

		Fix $x\notin\mathcal N$ and $0<r<R_*$. Choose $x_0\in D$ and a rational
		$R$ such that
		\[
		d(x_0,x)<r/16,
		\qquad r/2<R<3r/4.
		\]
		Then $x\in B(x_0,R/4)$ and $B(x_0,R)\subset B(x,r)$. A continuous path
		that crosses $B(x,r)$ before its lifetime must first cross $B(x_0,R)$.
		Set $K=2^\beta/a$, where $a$ is from
		Lemma~\ref{lem:fixed-ball}. Since $Kr^{-\beta}\ge a^{-1}R^{-\beta}$,
		we have for each path,
		\[
		e^{-Kr^{-\beta}S(x,r)}
		\le e^{-a^{-1}R^{-\beta}S(x_0,R)}.
		\]
		Taking expectations and using \eqref{eq:fixed-ball-laplace} proves \eqref{eq:uniform-laplace}.
	\end{proof}

	We can improve the above probabilistic estimate to a sub-Gaussian tail estimate.

	\begin{proposition}\label{prop:tail}
		There are constants $C_1,c_1>0$ such that
		\begin{equation}\label{eq:crossing-tail}
			\PP_x(S(x,R)\le t)
			\le C_1\exp\!\left[
			-c_1\left(\frac{R^\beta}{t}\right)^{1/(\beta-1)}
			\right]
		\end{equation}
		for every $x\notin\mathcal N$, every $R>0$, and every $0<t<r_0^\beta$.
	\end{proposition}

	\begin{proof}
		First suppose that
		\begin{equation}\label{eq:lambda-range}
			\lambda>\lambda_*:=KR_*^{-\beta},
		\end{equation}
		where $\lambda_*=0$ if $R_*=\infty$, and let
		$\ell=(K/\lambda)^{1/\beta}<R_*$.
		Fix $x\notin\mathcal N$ and define
		\[
		\tau_j=S(x,j\ell),\qquad j\ge1,\qquad \tau_0=0.
		\]
		If $\tau_j<\infty$, then $\tau_j<\zeta$ and
		$z=X_{\tau_j}\notin\mathcal N$. To reach the next radial sphere,
		a path must first reach distance at least $\ell$ from $z$ due to the triangle inequality
		\[
		d(z,w)\ge d(x,w)-d(x,z)\ge\ell
		\quad\text{when }d(x,w)\ge(j+1)\ell.
		\]
		The strong Markov property and Proposition~\ref{prop:uniform-laplace}
		therefore give
		\begin{equation}\label{eq:iteration}
			\EE_x[e^{-\lambda\tau_{j+1}}]
			\le\frac12\EE_x[e^{-\lambda\tau_j}].
		\end{equation}
		The assertion includes $j=0$.

		For $n=\lfloor R/\ell\rfloor$, continuity gives $\tau_n\le S(x,R)$
		on the event $\{S(x,R)<\infty\}$. Iterating
		\eqref{eq:iteration} gives
		\begin{equation}\label{eq:laplace-global}
			\EE_x[e^{-\lambda S(x,R)}]
			\le2^{-n}
			\le2\exp[-c_2R\lambda^{1/\beta}],
			\qquad c_2=(\log2)K^{-1/\beta}.
		\end{equation}
		For $n=0$, the same inequality holds by using
		$\EE_x[e^{-\lambda S(x,R)}]\le1$.

		Consequently, for $\lambda>\lambda_*$,
		\begin{equation}\label{eq:chernoff}
			\PP_x(S(x,R)\le t)
			\le2\exp[\lambda t-c_2R\lambda^{1/\beta}].
		\end{equation}
		The minimizer of the right-hand side is
		\[
		\lambda_{\mathrm{opt}}
		=\left(\frac{c_2R}{\beta t}\right)^{\beta/(\beta-1)}.
		\]
		If $\lambda_{\mathrm{opt}}>\lambda_*$, substitution in
		\eqref{eq:chernoff} gives \eqref{eq:crossing-tail} with
		\[
		c_1=(\beta-1)(c_2/\beta)^{\beta/(\beta-1)}.
		\]
		If $\lambda_{\mathrm{opt}}=\lambda_*$, apply \eqref{eq:chernoff} with $\lambda=\lambda_*+1/n$ and let $n\to\infty$. Continuity of $\lambda t-c_2R\lambda^{1/\beta}$ gives the same estimate.

		When $r_0=\infty$, all $R,t>0$ are covered. When $r_0<\infty$ and
		$\lambda_{\mathrm{opt}}\le\lambda_*$, we have
		\[
		\frac{R}{t^{1/\beta}}
		\le\frac{\beta}{c_2}(\lambda_*t)^{(\beta-1)/\beta}
		\le\frac{\beta}{c_2}
		(K\varepsilon^{-\beta})^{(\beta-1)/\beta},
		\]
		as $t<r_0^\beta$ and $R_*=\varepsilon r_0$.
		Thus $(R^\beta/t)^{1/(\beta-1)}$ is bounded by a constant. Increasing $C_1$, the trivial bound by $1$ proves
		\eqref{eq:crossing-tail} in this remaining case.
	\end{proof}

	We are now in a position to prove Theorem~\ref{thm:main1}.
	\begin{proof}[Proof of Theorem~\ref{thm:main1}]
		Fix $0<t<r_0^\beta$. For almost every pair $(x,y)$, we have
		\begin{equation}\label{eq:kernel-semigroup}
			p_t(x,y)=\int_Mp_{t/2}(x,z)p_{t/2}(y,z)\dd\mu(z).
		\end{equation}
		Let $D=d(x,y)$. If $D=0$, just apply \cite[Theorem 3.3]{GHL08}.
		Otherwise every $z\in M$ satisfies
		$d(x,z)\ge D/2$ or $d(y,z)\ge D/2$.
		Splitting \eqref{eq:kernel-semigroup} accordingly and applying
		\cite[Theorem 3.3]{GHL08} implies
		\begin{align}
			p_t(x,y)\le Ct^{-\alpha/\beta}\Bigl(&
			\int_{\{d(x,z)\ge D/2\}}p_{t/2}(x,z)\dd\mu(z)+\int_{\{d(y,z)\ge D/2\}}p_{t/2}(y,z)\dd\mu(z)
			\Bigr).
			\label{eq:split-kernel}
		\end{align}
		For almost every starting point, these integrals are the transition probabilities
		\[
		\PP_x\bigl(d(x,X_{t/2})\ge D/2,\ t/2<\zeta\bigr)
		\quad\text{and}\quad
		\PP_y\bigl(d(y,X_{t/2})\ge D/2,\ t/2<\zeta\bigr).
		\]
At the fixed time $s=t/2$, let
\[
Q_s(x,A)=\PP_x(X_s\in A,\ s<\zeta),\qquad
K_s(x,A)=\int_Ap_s(x,z)\,d\mu(z).
\]
Choose a countable algebra generating the Borel sets of $M$. Proper association, followed by bounded monotone approximation where needed, gives $Q_s(x,A)=K_s(x,A)$ for almost every $x$ for each member $A$ of this algebra. Outside the union of these countably many null sets, the two finite measures agree on the entire algebra and hence on every Borel set. The exceptional set is therefore independent of $A$. We may consequently take $A=\{z:d(x,z)\ge d(x,y)/2\}$, even though it depends on $(x,y)$.

		By continuity before the lifetime, the displayed events imply
		$S(x,D/2)\le t/2$ and $S(y,D/2)\le t/2$, respectively.
		Proposition~\ref{prop:tail} and \eqref{eq:split-kernel} now give
		\[
		p_t(x,y)\le Ct^{-\alpha/\beta}
		\exp\!\left[-c\left(\frac{D^\beta}{t}\right)^{1/(\beta-1)}\right].
		\]
		Since a properly exceptional set has measure zero, the exclusions
		made above do not change the asserted almost-everywhere conclusion.
		The proof is complete.
	\end{proof}

\section{Proof of Theorem \ref{thm:main2}}

For a space $X$, we use the superscript $X$ to indicate that an object is defined on
$X$. For example, $(\cE^X,\cF^X)$ denotes a Dirichlet form on $X$.

\begin{figure}[ht]
	\centering
		\subfloat[]{%
			\begin{tikzpicture}[inner sep=0pt]
				\node at (0,0) {\includegraphics[width=15pt]{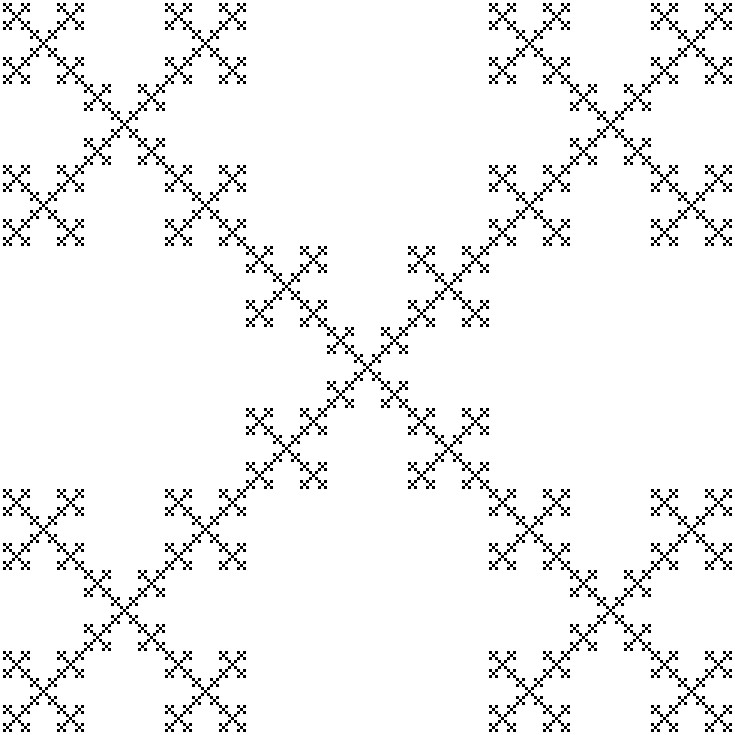}};
		\end{tikzpicture}}
	\hspace{5em}%
		\subfloat[]{
		\begin{tikzpicture}
			\foreach\w in {15pt}{
				\foreach\x in {0,1,2}{
					\foreach\y in {0,1,2}{
						\pgfmathparse{int(mod(\x+\y,2))}
						\let\r\pgfmathresult
						\ifnum\r=0
						\foreach\xx in {0,1,2}{
							\foreach\yy in {0,1,2}{
								\pgfmathparse{int(mod(\xx+\yy,2))}
								\let\rr\pgfmathresult
								\ifnum\rr=0
								\node at (\x*3*\w+\xx*\w,\y*3*\w+\yy*\w) {\includegraphics[width=\w]{vc}};
								\fi
							}
						}
						\fi
					}
				}
			}
		\end{tikzpicture}}	\hspace{5em}
	\subfloat[]{%
		\begin{tikzpicture}[inner sep=0pt]
			\def\w{15pt}                     
			\foreach \i in {0,1,2,3,4,5,6,7,8} {%
				\node at (\i*\w,\i*\w) {\includegraphics[width=\w]{vc}};
			}
	\end{tikzpicture}}
	\caption{(a) The Vicsek set K, (b) The unbounded Vicsek set Y, (c) The union of all Vicsek cells in $Y$ that intersect the diagonal ray, that is, Z} \label{vc}
\end{figure}

The construction is quite lengthy. Recall that the Vicsek set is the unique non-empty compact subset $K\subseteq\R^2$ such that
$K=\cup_{i=0}^4F_i(K)$ (See Figure \ref{vc}(a)), where \[
F_i(x)=\frac{1}{3}(x-p_i)+p_i,\qquad x\in\R^2,\ i=0,\dots,4
\]
with
\[
p_0=(0,0),\quad p_1=(1,0),\quad p_2=(1,1),\quad p_3=(0,1),\quad
p_4=\left(\tfrac{1}{2},\tfrac{1}{2}\right).
\]

For any $n\ge 0$, let
\[
W_n=\{w=w_1\cdots w_n: w_i\in\{0,\dots,4\} \text{ for all } i=1,\dots,n\}.
\]
For any $w=w_1\cdots w_n\in W_n$, let $|w|=n$ be the length of
the word $w$, and let
\[
F_w=F_{w_1}\circ\cdots\circ F_{w_n},\qquad K_w=F_w(K).
\]
Let $m^K$ be the unique positive Radon measure on $K$ such that
$m^K(K_w)=5^{-|w|}$ for any $w\in \cup_{n\ge 0}W_n$. Define $V_0=\{p_0,\dots,p_4\}$ and
$V_{n+1}=\cup_{i=0}^4F_i(V_n)$ for $n\ge 0$, then $V_*:=\cup_{n\ge 0}V_n$ is dense in
$K$. For any $n\ge 0$ and any $u\in l(V_n)$, define
\[
\cE^K_n(u,u)=3^n\sum_{w\in W_n}\sum_{i=0}^{3}
\left(u(F_w(p_i))-u(F_w(p_4))\right)^2.
\]
Then, for any $u\in l(V_*)$, the sequence $\{\cE^K_n(u,u)\}_{n\ge 0}$ is non-decreasing. Let
\[
\cE^K(u,u)=\lim_{n\to+\infty}\cE^K_n(u,u)=\sup_{n\ge 0}\cE^K_n(u,u),
\]
\[
\cF^K=\left\{u\in C(K): \sup_{n\ge 0}\cE^K_n(u,u)<+\infty\right\}.
\]
Then $(\cE^K,\cF^K)$ is a strongly local regular Dirichlet form on $L^2(K;m^K)$.

For $a\in\R$, $y\in\R^2$, and $A\subseteq\R^2$, write $aA+y:=\{ax+y: x\in A\}$. We call $A$ a Vicsek cell if $A=K+y$ for some $y\in\R^2$. For a Vicsek cell $A=K+y$, define $m^A$ and $(\cE^A,\cF^A)$ as the
pushforwards under the map $x\mapsto x+y$ of $m^K$ and $(\cE^K,\cF^K)$, respectively.
Let
\[
Y:=\bigcup_{n\ge 0}3^nK
\]
be the unbounded Vicsek set (See Figure \ref{vc}(b)), then $Y$ is a union of countably many Vicsek cells, and the
intersection of any two distinct Vicsek cells $A_1,A_2\subseteq Y$ is either empty or a
singleton. Let $d^Y$ be the geodesic metric on $Y$, then $d^Y$ is bi-Lipschitz
equivalent to the Euclidean metric. Let $m^Y$ be the unique positive Radon measure on
$Y$ such that $m^Y|_A=m^A$ for any Vicsek cell $A\subseteq Y$. Then $(Y,\text{ }d^Y,\text{ }m^Y)$ is an
unbounded complete geodesic metric measure space satisfying $\VC(\alpha)$.

Define
\[
\cE^Y(u,u)=\sum_{\substack{A\subseteq Y\\ A \text{ is a Vicsek cell}}}
\cE^A(u|_A,u|_A),
\]
\[
\cF^Y=\left\{u\in C(Y)\cap L^2(Y,m^Y): u|_A\in\cF^A \text{ for any Vicsek cell } A\subseteq Y,\
\cE^Y(u,u)<+\infty\right\}.
\]
Then $(\cE^Y,\cF^Y)$ is a strongly local regular Dirichlet form on $L^2(Y;m^Y)$.

Let
\[
Z:=\bigcup_{n\ge 0}\left(K+(n,n)\right)
\]
be the union of all Vicsek cells in $Y$ that intersect the diagonal ray
$\{(x,x): x\in[0,+\infty)\}\subseteq\R^2$ (See Figure \ref{vc}(c)), which is a connected closed subset of
$(Y,d^Y)$. Let $d^Z=d^Y|_{Z\times Z}$; then $d^Z$ is the geodesic metric on $Z$. Let $m^Z=m^Y|_Z$. We equip $Z$ with the natural tree form $(\cE^Z,\cF^Z)$, defined by summing the energies of its unit cells, with matching values at their common vertices, exactly as in the definition of $(\cE^Y,\cF^Y)$.
Then $(Z,d^Z,m^Z)$ is an unbounded complete geodesic metric measure space satisfying
$\VC(\Phi^Z)$ with
\[
\Phi^Z(r)=
\begin{cases}
r^\alpha, & \text{if } r\in(0,1],\\
r, & \text{if } r\in(1,+\infty).
\end{cases}
\]

 The volume estimates and the geodesic property, together with standard heat kernel results, show that $\HK(\beta)$ holds on $Y$ and that
$\HK(\Psi^Z)$ holds on $Z$, where
\[
\Psi^Z(r):=
\begin{cases}
r^\beta, & \text{if } r\in(0,1],\\
r^2, & \text{if } r\in(1,+\infty).
\end{cases}
\]See, for example, \cite[Theorem~3.1 and Corollary~3.1]{Kum}.

Let $S=Y\times Y$ be equipped with the metric
\[
d^S\left((x_1,x_2),(y_1,y_2)\right):=\max\left\{d^Y(x_1,y_1),d^Y(x_2,y_2)\right\}
\]
and the product measure $m^S:=m^Y\otimes m^Y$, then $(S,d^S,m^S)$ is an unbounded complete geodesic
metric measure space satisfying $\VC(2\alpha)$. The space $S$ admits the following product heat kernel
\[
 p_t^{Y\times Y}((x_1,x_2),(y_1,y_2))=p_t^Y(x_1,y_1)p_t^Y(x_2,y_2),
\]
where $p_t^Y$ is the heat kernel on $Y$; the product kernel satisfies $\HK(\beta)$ by a direct calculation. The corresponding $(\cE^S,\cF^S)$ is a strongly local regular Dirichlet
form on $L^2(S;m^S)$. Following the arguments in \cite[Theorem~14.20 and Example~14.12]{G} and using $\VC(2\alpha)$ and $\HK(\beta)$, there exists
$C_F>0$ such that for any non-empty open subset $\Omega\subseteq S$,
\begin{equation}\label{slow-spec}
\lambda^S_{\min}(\Omega)\ge\frac{C_F}{m^S(\Omega)^{\frac{\beta}{2\alpha}}}.
\end{equation}
Similarly, let $F=Z\times Z$ be equipped with the metric
\[
d^F\left((x_1,x_2),(y_1,y_2)\right):=\max\left\{d^Z(x_1,y_1),d^Z(x_2,y_2)\right\}
\]
and the product measure $m^F:=m^Z\otimes m^Z$.  Then $(F,d^F,m^F)$ is an unbounded complete geodesic
metric measure space satisfying $\VC(\Phi^F)$, where
\begin{equation}\label{volume-fast}
\Phi^F(r)=\Phi^Z(r)^2=
\begin{cases}
r^{2\alpha}, & \text{if } r\in(0,1],\\
r^2, & \text{if } r\in(1,+\infty).
\end{cases}
\end{equation}
Similarly, it admits a heat kernel associated with a strongly local regular Dirichlet form $(\cE^F,\cF^F)$  on $L^2(F;m^F)$
satisfying $\HK(\Psi^F)$, where
\begin{equation}\label{time-fast}
\Psi^F(r)=\Psi^Z(r)=
\begin{cases}
r^\beta, & \text{if } r\in(0,1],\\
r^2, & \text{if } r\in(1,+\infty).
\end{cases}
\end{equation}
Moreover, by $\VC(\Phi^F)$ and $\HK(\Psi^F)$, there exists $C_F>0$ such that for any
non-empty open subset $\Omega\subseteq F$,
\begin{equation}\label{fast-spec}
\lambda^F_{\min}(\Omega)\ge
\begin{cases}
\frac{C_F}{m^F(\Omega)^{\frac{\beta}{2\alpha}}}, & \text{if } m^F(\Omega)<1,\\[6pt]
\frac{C_F}{m^F(\Omega)}, & \text{if } m^F(\Omega)\ge 1.
\end{cases}
\end{equation}

We use the letters $S$ and $F$ for the `slow' and `fast' product-fractals, respectively.
Let $D_0=\{(n,n):n\in\mathbb N\}$ and $D=D_0\times D_0\subset F\subset S$.
We use tagged copies of $S$ and $F$ to form the disjoint union
\[
\widetilde X:=(S\times\{\mathrm s\})\cup(F\times\{\mathrm f\}),
\]
where $\mathrm s$ and $\mathrm f$ are distinct labels. Define an equivalence relation on $\widetilde X$ by
\[
(x,i)\sim(y,j)
\quad\Longleftrightarrow\quad
x=y\ \text{and}\ \bigl(i=j\ \text{or}\ x\in D\bigr).
\]
The glued space is the quotient
\[
X:=\widetilde X/\!\sim.
\]
Thus $(z,\mathrm s)$ and $(z,\mathrm f)$ represent the same point of $X$ precisely when $z\in D$; the two copies of every $z\in F\setminus D$ remain distinct.

Define the canonical injective maps
\[
\pi^S:S\longrightarrow X,\quad x\longmapsto[(x,\mathrm s)],
\qquad
\pi^F:F\longrightarrow X,\quad y\longmapsto[(y,\mathrm f)],
\]
where $[(x,i)]$ denotes the equivalence class of $(x,i)$.

From this point onwards, we identify each sheet with its image in $X$.
The two images intersect precisely in $\pi^S(D)=\pi^F(D)$, which we also denote by $D$; under this convention, $S\cap F=D$ in $X$.

For points $x,y$ in the same sheet $S$ (or $F$), we simply define the metric
\[
d^X(x,y):=d^S(x,y)\ (\text{or}\ d^F(x,y)).
\]
For points in different sheets ($x\in S$, $y\in F$), the quotient distance is defined as
\[
d^X(x,y):=\min_{z\in D}\{d^S(x,z)+d^F(z,y)\}.
\]

Intuitively, $X$ is made up by gluing $S$ and $F$ along countably many points in $D$. It is easy to see that $X$ is geodesic and complete.

\begin{figure}[htbp]\centering
\begin{tikzpicture}[x=1cm,y=1cm]
\draw[very thick,blue!65!black,domain=0:8,samples=101] plot (\x,{0.6*(1-cos(90*\x))});
\draw[very thick,orange!75!black] (0,0)--(2,-1)--(4,0)--(6,-1)--(8,0);
\foreach \x in {0,4,8}{\fill (\x,0) circle (2pt);}
\node at (2,1.5){$S$};\node at (2,-1.3){$F$};\node at (4,-.4){$D$};
\end{tikzpicture}
\caption{Schematic illustration of the quotient space $X$}
\end{figure}
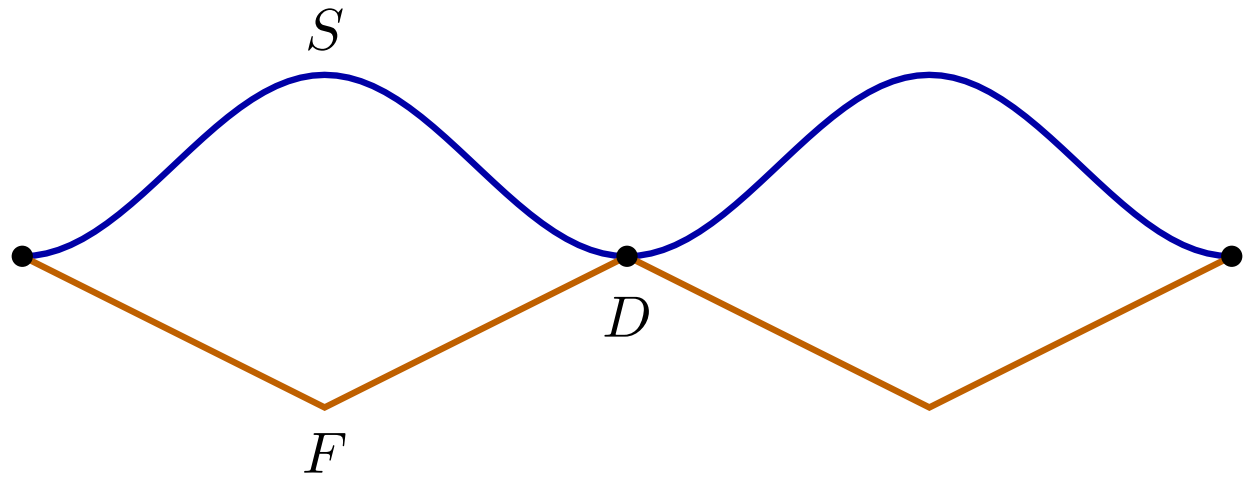
Since $2\alpha>\beta$, it follows from $\VC(2\alpha)$, $\HK(\beta)$, $\VC(\Phi^F)$ and
$\HK(\Psi^F)$ that\begin{equation}\label{polar}
\Capm^S_1(D)=m^S(D)=\Capm^F_1(D)=m^F(D)=0.
\end{equation}

By \cite[Theorem~1.2]{GHL15}, the two-sided heat kernel bounds give the capacity upper bounds on each sheet. We have the following.
\begin{proposition}
Conditions $\Capm(\beta)$ and $\Capm(\Psi^F)$ are satisfied on $S$ and $F$, respectively.
\end{proposition}

For $W=S,F$, a cutoff in a ball of radius $r<1$ has energy at most $Cr^{2\alpha-\beta}$ and squared $L^2$ norm at most $Cr^{2\alpha}$. Thus $\Capm_1^W(\{z\})=0$ for every $z\in D$, and countable subadditivity gives \eqref{polar}.
For any subset $A\subseteq X$, let
\[
A^S:=(\pi^S)^{-1}(A),\qquad A^F:=(\pi^F)^{-1}(A).
\]
Under the identification introduced above, we also write $A^S=A\cap S$ and $A^F=A\cap F$.

Let
\[
m^X=(\pi^S)_{\#}m^S+(\pi^F)_{\#}m^F.
\]
It follows that $m^X$ is a positive Radon measure on $X$, and for any measurable subset
$A\subseteq X$, $m^X(A)=m^S(A^S)+m^F(A^F)$, with
$L^2(X;m^X)=L^2(S;m^S)\oplus L^2(F;m^F)$.
For any measurable function $u$ on $X$, let
\[
u^S:=u\circ\pi^S,\qquad u^F:=u\circ\pi^F.
\]
Under the same identification, we also denote these functions by $u|_S$ and $u|_F$, respectively.

Let
\[
\cE^X(u,u)=\cE^S(u^S,u^S)+\cE^F(u^F,u^F),\qquad
\cF^X=\cF^S\oplus\cF^F,
\]
then it is easy to see that $(\cE^X,\cF^X)$ is a strongly local regular Dirichlet form on $L^2(X;m^X)$ as the points in $D$ have zero $1$-capacity. For
any $u\in\cF^X$ and any Borel measurable set $A\subseteq X$, the corresponding energy measure is given by
\[
\Gamma^X(u)(A)=\Gamma^S(u^S)(A^S)+\Gamma^F(u^F)(A^F).
\]

\begin{lemma}
The metric measure space $(X,d^X,m^X)$ satisfies $\VC(2\alpha)$.
\end{lemma}

\begin{proof}
Let $B^X(x,r)$ be a ball in $X$. Recall $\VC(2\alpha)$ on $S$ and $\VC(\Phi^F)$ on $F$. When $x\in S$, by using the fact that
$B^X(x,r)\supseteq B^S(x,r)$, we have
\[
m^X\left(B^X(x,r)\right)\ge m^S\left(B^S(x,r)\right)\asymp r^{2\alpha}.
\]
If $B^X(x,r)\cap D=\emptyset$, then $B^X(x,r)=B^S(x,r)$; otherwise, choose
$y\in B^X(x,r)\cap D$, then
$B^X(x,r)\subseteq B^S(x,r)\cup B^F(y,2r)$. The desired bound is immediate in the first case; in the second case, we have
\[
m^X\left(B^X(x,r)\right)\le m^S\left(B^S(x,r)\right)+m^F\left(B^F(y,2r)\right)
\asymp r^{2\alpha}+\Phi^F(2r)\asymp r^{2\alpha}.
\]
Every point of $F$ is at distance at most some fixed $L\ge1$ from $D$, since every unit cell has bounded diameter. When $x\in F$, we distinguish two cases for the lower bound. If $r<2L$, then $B^X(x,r)\supseteq B^F(x,r)$, and so
\[
m^X\left(B^X(x,r)\right)\ge m^F\left(B^F(x,r)\right)\asymp\Phi^F(r)\asymp r^{2\alpha}.
\]
If $r\ge 2L$, choose $y\in D$ with
$d^X(x,y)=d^F(x,y)\le L$. Then $B^X(x,r)\supseteq B^S(y,r-L)$, and so
\[
m^X\left(B^X(x,r)\right)\ge m^S\left(B^S\left(y,r-L\right)\right)
\asymp r^{2\alpha}.
\]
For the upper bound, if $B^X(x,r)\cap D=\emptyset$, then $B^X(x,r)=B^F(x,r)$;
otherwise, choose $y\in B^X(x,r)\cap D$, then
$B^X(x,r)\subseteq B^F(x,r)\cup B^S(y,2r)$. The desired bound is immediate in the first case; in the second case,
\[
m^X\left(B^X(x,r)\right)\le m^F\left(B^F(x,r)\right)+m^S\left(B^S(y,2r)\right)
\asymp\Phi^F(r)+(2r)^{2\alpha}\asymp r^{2\alpha}.
\]
\end{proof}

\begin{lemma}
The capacity upper bound $\Capm(\beta)_{\leq}$ holds for $(\cE^X,\cF^X)$.
\end{lemma}

\begin{proof}
Let $A_1,A_2>1$ be the constants in $\Capm(\beta)_{\leq}$ on $S$,
$\Capm(\Psi^F)_{\leq}$ on $F$, respectively. Let $B^X(x,r)$ be a ball in $X$. When
$x\in S$, by $\Capm(\beta)_{\leq}$ and $\VC(2\alpha)$ on $S$, there exists a cutoff
function $\phi^S\in\cF^S$ for $B^S(x,r)\subseteq B^S(x,A_1r)$ such that
\[
\cE^S(\phi^S,\phi^S)\lesssim\frac{m^S\left(B^S(x,r)\right)}{r^\beta}
\asymp r^{2\alpha-\beta}.
\]
If $B^X(x,r)\cap D=\emptyset$, then let $\phi^F=0\in\cF^F$; otherwise, choose
$y\in B^X(x,r)\cap D$, then $B^X(x,r)\subseteq B^S(x,r)\cup B^F(y,2r)$, by
$\Capm(\Psi^F)_{\leq}$ and $\VC(\Phi^F)$ on $F$, there exists a cutoff function
$\phi^F\in\cF^F$ for $B^F(y,2r)\subseteq B^F(y,2A_2r)$ such that
\[
\cE^F(\phi^F,\phi^F)\lesssim\frac{m^F\left(B^F(y,2r)\right)}{\Psi^F(2r)}
\asymp\frac{\Phi^F(r)}{\Psi^F(r)}
=
\begin{cases}
r^{2\alpha-\beta}, & \text{if } r<1,\\
1, & \text{if } r\ge 1.
\end{cases}
\]
Let $\phi$ be defined on $X$ according to $\phi^S$ on $S$ and $\phi^F$ on $F$, where the values on the polar set $D$ are immaterial. Then $\phi\in\cF^X$, $\widetilde\phi=1$ quasi-everywhere on $B^X(x,r)$, its support is contained in $B^X(x,4(A_1+A_2)r)$, and
\[
\cE^X(\phi,\phi)=\cE^S(\phi^S,\phi^S)+\cE^F(\phi^F,\phi^F)
\lesssim r^{2\alpha-\beta}+\frac{\Phi^F(r)}{\Psi^F(r)}
\asymp r^{2\alpha-\beta}
\asymp \frac{m^X\left(B^X(x,r)\right)}{r^\beta},
\]
where we use the fact that $2\alpha>\beta$ and
 $\VC(2\alpha)$ on $X$. Similarly, for $x\in F$, the same construction $\phi$ satisfies
 \[
\cE^X(\phi,\phi)
\lesssim r^{2\alpha-\beta}+\Phi^F(2r)/\Psi^F(2r)\lesssim m^X(B^X(x,r))/r^\beta.
\] Applying the same construction with inner radius $2r$ gives a cutoff equal to $1$ quasi-everywhere on $\overline{B^X(x,r)}$, with the same energy bound up to a constant. Therefore $\Capm(\beta)_{\leq}$ holds for $(\cE^X,\cF^X)$.
\end{proof}

\begin{lemma}
The Faber--Krahn inequality $\FK(\beta)$ holds for $(\cE^X,\cF^X)$.
\end{lemma}

\begin{proof}
We prove $\FK(\beta)$ holds for $(\cE^X,\cF^X)$ with $\nu=\frac{1}{2}$. Let
$B^X(x,r)$ be a ball in $X$ and let $\Omega\subseteq B^X(x,r)$ be a non-empty open
subset. Since $\beta=\alpha+1$, by $\VC(2\alpha)$ on $X$, it suffices to show that for
any $u\in\cF^X(\Omega)$,
\[
\|u\|^2_{L^2(X;m^X)}\lesssim rm^X(\Omega)^{\frac{1}{2}}\cE^X(u,u).
\]
Recall that $\Omega^S=\Omega\cap S$, $\Omega^F=\Omega\cap F$, $u^S=u|_S$, and
$u^F=u|_F$. Noting that $\Omega^S,\Omega^F$ are open subsets of $S,F$, respectively,
by \eqref{polar}, $u^S\in\cF^S(\Omega^S)$, $u^F\in\cF^F(\Omega^F)$. Moreover,
\[
m^X(\Omega)=m^S(\Omega^S)+m^F(\Omega^F),\qquad
\cE^X(u,u)=\cE^S(u^S,u^S)+\cE^F(u^F,u^F).
\]
We first show \begin{equation}\label{ym25}
\|u^S\|^2_{L^2(S;m^S)}\lesssim rm^X(\Omega)^{\frac{1}{2}}\cE^X(u,u).
\end{equation}
This is trivial when $\Omega^S=\emptyset$ as $u^S=0$, thus it remains to consider the case $\Omega^S\neq\emptyset$. By \eqref{slow-spec},
\[
\|u^S\|^2_{L^2(S;m^S)}\lesssim m^S(\Omega^S)^{\frac{\beta}{2\alpha}}
\cE^S(u^S,u^S).
\]
Take $y\in\Omega^S$, then $\Omega^S\subseteq B^S(y,2r)$, by $\VC(2\alpha)$ on $S$,
\[
m^S(\Omega^S)^{\frac{\beta}{2\alpha}}
=m^S(\Omega^S)^{\frac{\beta}{2\alpha}-\frac{1}{2}}m^S(\Omega^S)^{\frac{1}{2}}
\le m^S\left(B^S(y,2r)\right)^{\frac{\beta}{2\alpha}-\frac{1}{2}}m^X(\Omega)^{\frac{1}{2}}
\lesssim rm^X(\Omega)^{\frac{1}{2}},
\]
proving the assertion.

We then show \begin{equation}
    \label{ym26}\|u^F\|^2_{L^2(F;m^F)}\lesssim rm^X(\Omega)^{\frac{1}{2}}\cE^X(u,u).
\end{equation}
Again, this is trivial when $\Omega^F=\emptyset$ as $u^F=0$, and it remains to consider the case $\Omega^F\neq\emptyset$. Fix $z\in\Omega^F$, then $\Omega^F\subseteq B^F(z,2r)$.
By \eqref{fast-spec},
\begin{equation}
\|u^F\|^2_{L^2(F;m^F)}\lesssim
\begin{cases}
m^F(\Omega^F)^{\frac{\beta}{2\alpha}}\cE^F(u^F,u^F),
& \text{if } m^F(\Omega^F)<1,\\[4pt]
m^F(\Omega^F)\cE^F(u^F,u^F),
& \text{if } m^F(\Omega^F)\ge 1.
\end{cases}
\end{equation}
We distinguish three cases.
If
$m^F(\Omega^F)<1$ and $r<1$, then by $\VC(\Phi^F)$ on $F$,
\[
m^F(\Omega^F)^{\frac{\beta}{2\alpha}}
=m^F(\Omega^F)^{\frac{\beta}{2\alpha}-\frac{1}{2}}m^F(\Omega^F)^{\frac{1}{2}}
\le m^F\left(B^F(z,2r)\right)^{\frac{\beta}{2\alpha}-\frac{1}{2}}m^X(\Omega)^{\frac{1}{2}}
\lesssim rm^X(\Omega)^{\frac{1}{2}}.
\]
If $m^F(\Omega^F)<1$ and $r\ge 1$, then
\[
m^F(\Omega^F)^{\frac{\beta}{2\alpha}}
=m^F(\Omega^F)^{\frac{\beta}{2\alpha}-\frac{1}{2}}m^F(\Omega^F)^{\frac{1}{2}}
\le m^F(\Omega^F)^{\frac{1}{2}}\le m^X(\Omega)^{\frac{1}{2}}
\le rm^X(\Omega)^{\frac{1}{2}}.
\]
If $m^F(\Omega^F)\ge 1$, then $m^F(B^F(z,2r))\ge 1$ and so $r\gtrsim 1$ by using $\VC(\Phi^F)$ on $F$, which implies that
\[
m^F(\Omega^F)=m^F(\Omega^F)^{\frac{1}{2}}m^F(\Omega^F)^{\frac{1}{2}}
\le m^F\left(B^F(z,2r)\right)^{\frac{1}{2}}m^X(\Omega)^{\frac{1}{2}}
\lesssim rm^X(\Omega)^{\frac{1}{2}}.
\]
Combining \eqref{ym25} and \eqref{ym26}, we obtain
\[
\|u\|^2_{L^2(X;m^X)}=\|u^S\|^2_{L^2(S;m^S)}+\|u^F\|^2_{L^2(F;m^F)}
\lesssim rm^X(\Omega)^{\frac{1}{2}}\cE^X(u,u).
\]
\end{proof}

\begin{lemma}
	The generalized capacity condition
	$\mathrm{Gcap}_{\leq}$ fails for
	$(\cE^X,\cF^X)$.
\end{lemma}

\begin{proof}
	Suppose, to the contrary, that
	$\mathrm{Gcap}_{\leq}$ holds.
	Fix $A>1$ such that $A^2c_{V,F}\geq2C_{V,F}$,
	and let $o=(0,0,0,0)\in D\subseteq X$.
	
	For $r>1$, regularity gives a function
	$f^F\in\cF^F\cap C_c(F)$ such that
	$0\leq f^F\leq1$ and $f^F=1$ on an open
	neighborhood of $\overline{B^F(o,Ar)}$.
	Define $f$ on $X$ by $f=0$ on $S\setminus D$
	and $f=f^F$ on $F\setminus D$.
	The values on the polar set $D$ are immaterial.
	By \eqref{polar},
	$f\in\cF^X\cap L^\infty(X)$.
	
	Since $A$ is fixed,
	$\mathrm{Gcap}_{\leq}$ provides a cutoff
	function $\phi_r\in\cF^X$ for
	$B^X(o,r)\subseteq B^X(o,Ar)$ such that
	\begin{align}
		\int_{B^X(o,Ar)\setminus B^X(o,r)}
		|\tilde f|^2\,d\Gamma^X(\phi_r,\phi_r)
		&\lesssim
		\int_{B^X(o,Ar)\setminus B^X(o,r)}
		d\Gamma^X(f,f) \notag\\
		&\quad+
		\frac1{r^\beta}
		\int_{B^X(o,Ar)\setminus B^X(o,r)}
		|f|^2\,dm^X.
		\label{cs-test}
	\end{align}
	The implicit constant may depend on $A$,
	but is independent of $r$ and $f$.
	
	Let $\phi_r^F=\phi_r|_F$.
	Then $\phi_r^F\in\cF^F$ is a cutoff function
	for $B^F(o,r)\subseteq B^F(o,Ar)$.
\begin{align*}
	\int_{B^X(o,Ar)\setminus B^X(o,r)}
	|\tilde f|^2\,d\Gamma^X(\phi_r,\phi_r)
	&=
	\int_{B^F(o,Ar)\setminus B^F(o,r)}
	|f^F|^2\,d\Gamma^F(\phi_r^F,\phi_r^F)\\
	&=\cE^F(\phi_r^F,\phi_r^F).
\end{align*}
Furthermore,
\begin{align*}
	\cE^F(\phi_r^F,\phi_r^F)
	&\geq
	\Capm^F\bigl(B^F(o,r),
	F\setminus B^F(o,Ar)\bigr)\\
	&\geq
	c\frac{m^F(B^F(o,r))}{\Psi^F(r)}
	\asymp1,
\end{align*}
where we used $\Capm(\Psi^F)_{\geq}$.

	On the other hand, the choice of $f$ and
	strong locality imply
	\[
	\int_{B^X(o,Ar)\setminus B^X(o,r)}
	d\Gamma^X(f,f)=0.
	\]

	Moreover, since $A^2c_{V,F}\geq2C_{V,F}$,
	\begin{align*}
		\tfrac12c_{V,F}A^2r^2
		&\leq(c_{V,F}A^2-C_{V,F})r^2
		\leq m^F(B^F(o,Ar)\setminus B^F(o,r))
		\leq C_{V,F}A^2r^2,\\
		\frac1{r^\beta}
		\int_{B^X(o,Ar)\setminus B^X(o,r)}|f|^2\,dm^X
		&=\frac{m^F(B^F(o,Ar)\setminus B^F(o,r))}{r^\beta}
		\asymp r^{2-\beta}.
	\end{align*}
	Therefore, \eqref{cs-test} yields
	\[
	1\lesssim r^{2-\beta}.
	\]
	Since $\beta>2$, letting $r\to+\infty$
	gives a contradiction.
\end{proof}

Theorem~\ref{thm:main2} follows immediately from the above lemmas.

The construction disproves Conjecture 2. Indeed, the semigroup on $X$ is the direct {sum} of those on $S$ and $F$, so its heat kernel agrees with $p_t^F$ on $F\times F$ and is zero between the two sheets (up to null sets). By $\HK(\Psi^F)$, there are $c,\eta>0$ such that, for $t\ge1$,
\[
p_t^X(x,y)=p_t^F(x,y)\ge c t^{-1}
\quad\text{for almost every }x,y\in F\text{ with }d^F(x,y)<\eta\sqrt t.
\]
These pairs have positive product measure. On the other hand, $\VC(2\alpha)$ on $X$ and a sub-Gaussian upper estimate with walk-dimension $\beta$ would give
\[
p_t^X(x,y)\le C t^{-2\alpha/\beta}.
\]
Since $2\alpha/\beta>1$, these inequalities contradict each other as $t\to\infty$.
\section*{Acknowledgement}
The authors thank Aobo Chen and Meng Yang for providing great ideas and examples.
	Jin Gao was supported by the National Natural Science Foundation of China (Grant No. 12271282) and the Zhejiang Provincial Natural Science Foundation of China (Grant No. LQN25A010019). Zhenyu Yu was supported by the Natural Science Foundation of Hunan Province, China (Grant No. 2025JJ60039) and the National University of Defense Technology (Grant No. ZK25-05).
\section*{Use of AI}
The authors used OpenAI's GPT-5.6 and GPT-6.0 to assist with preliminary mathematical exploration, language polishing, and manuscript formatting. The authors take full responsibility for the entire content of the manuscript, including the correctness of the mathematical arguments and results and the accuracy of the references.

\end{document}